\documentclass[11pt,a4paper]{amsart}

\allowdisplaybreaks
\usepackage{amsmath,amssymb,amsthm,graphicx,hyperref,color}
\usepackage{booktabs}
\usepackage{multirow}
\usepackage{placeins}
\usepackage{float}
\newtheorem{remark}{Remark}
\hypersetup{colorlinks=true,
	linkcolor=blue,
	filecolor=magenta,
	urlcolor=cyan,}

\title[PINN and open sphere visualization]{A Physics-Informed Neural Network for Geodesic Computation and Visualization of Geodesic Spheres on Riemannian Manifolds}

\author[S. Bhattacharyya]{Sujit Bhattacharyya}
\begin{document}
	
	\maketitle
	
	\begin{abstract}
		Visualizing geodesic structures associated with general Riemannian metrics is a challenging problem in computational differential geometry. In this work, we present a computational framework for visualizing geodesics and geodesic spheres on a three-dimensional Riemannian manifold equipped with the metric $g=\text{diag}(1,x^2,1), x>0$. The associated Christoffel symbols and curvature tensors are computed symbolically, and the resulting geodesic equations are solved using both \texttt{SciPy}'s numerical integration method and a Physics-Informed Neural Network (PINN). The two solutions are compared to assess the accuracy of the neural network approximation. We further develop geodesic spike diagrams and surface visualizations to illustrate the deformation of geodesic spheres relative to their Euclidean counterparts. Intrinsic Riemannian surface area and volume are estimated numerically using the induced metric and the Riemannian volume density, respectively. For the chosen metric and a geodesic ball of radius one, the computed area and volume are approximately 
		$0.362\%$ and $0.364\%$ smaller than the corresponding Euclidean values at the base grid resolution. A grid-refinement study is also performed to examine the sensitivity of these estimates to numerical resolution. The proposed framework demonstrates how symbolic differential geometry, numerical integration, and PINNs can be combined to investigate and visualize geometric structures associated with Riemannian metrics.
	\end{abstract}
	\tableofcontents
	\section{Introduction}
	\subsection{Statement of the problem}
		Imagination comes when one visualize something and then extends it in mind. Visualization plays a very important role in mathematics and specially in geometry. Objects that appear in Euclidean geometry are easy to visualize. Any function $f:\mathbb{R}\to\mathbb{R}$ can be plotted as a subset $G(f)=\{(x,f(x)):x\in\mathbb{R}\}\subset \mathbb{R}^2$. In similar fashion a map $f:\mathbb{R}^2\to \mathbb{R}$ can be visualized using the graph $G(f)=\{(x,y,f(x,y)):(x,y)\in\mathbb{R}^2\}\subset \mathbb{R}^3$. Things gets interesting when dimension reaches to $4$. For example, a map $f:\mathbb{R}^3\to \mathbb{R}$ or $f:\mathbb{R}^2\to\mathbb{R}^2$ can produce a graph $G(f)$ as a subset of $\mathbb{R}^4$. We have no direct access to the 4 dimensional space and hence they can not be visualized in the traditional manner. Another obstacle is the metric defined on the space. For example, we are able to plot those graphs on $\mathbb{R}^2$ or $\mathbb{R}^3$ under the Euclidean metric. Now what if the metric changes? Can we visualize them on the space again? One possible answer is that if $(\mathbb{R},d_{u})$ is the usual metric space then the space $(\mathbb{R},\lambda d_u)$ is the scaled metric space whose open unit spheres are going to be $\lambda$ times bigger/smaller than the original space. But this explanation is going to work simple looking metrics. \textbf{We are talking about visualizing open unit spheres of Riemannian manifolds.} One may say that yes we can see it; the sphere $\mathbb{S}^2$, the cube $[0,1]^3$, we can draw them and it has its own Riemannian metric. But my question is \textbf{what if we have a Riemannian metric, for example $ds^2 = x^2 dx^2+\sin y dy^2$, how to ``see" the open unit sphere under this metric?}
		
		In this article we are going to visualize the open unit spheres of a Riemannian manifold equipped with an abstract\footnote{The word abstract means it is not a classical metric. For example, $g$ is not hyperbolic, spherical or Euclidean metric. The functions $g_{ij}$ are customizable.} Riemannian metric $g$ by solving the geodesic equations. We are also going to explore a volume comparison between the open sphere of the manifold and the open Euclidean sphere. The geodesic equations that are arising from the metric can be solved via Python-\texttt{SciPy}'s integrate IVP module, but it will produce a discrete solution, however if we configure \texttt{SciPy} then it will interpolate between the points and provide a continuous solution. But we create a Physics-Informed Neural Network (PINN) model which will solve the equations in a continuous manner. We are applying Neural Network to understand the nature of the solution in a continuous manner rather than having an interpolation of noisy data. This implementation helps us to study real world data and it also shows that PINN can be applied to solve geodesic equations over Riemannian manifolds. We also compare both the solution (PINN vs \texttt{SciPy}) to verify the accuracy and stability of the model.
		
		\subsection{Related works}
		This kind of work is not new in computational Riemannian geometry. In the past these works were done manually i.e., by hand. One can see the work of Hui and his coauthors \cite{Hui2011,HuiMatsuyama2015,Hui2016,HuiMatsuyama2016} for such investigations. But after the involvement of mathematical softwares like Mathematica, Matlab and Python especially after the introduction of symbolic computations several researchers have been investigating Riemannian manifolds and its various curvature tensors. One can see the work of Shaikh et al. \cite{ShaikhHuiDattaSarkar2024, ShaikhHuiSarkarBabu2024} for analysis on space-time metrics. For computations over vector field related equations we refer to the works of Blaga \cite{Blaga2025a, Blaga2025b, Blaga2025c}.On the other hand several nonlinear differential equations were being solved using PINN, in this regard the works of Sur \cite{SurLotfyBhattacharyyaDuttaDas2026, SurHussainCraciunSinghalRanjitYaylaciDas2026} is worth mentioning. Other research leads to theoretical understanding of the space where as this article provides a direct visualization of the unseen open spheres of the Riemannian space. To the best of the authors knowledge and belief this kind of study is new to computational Riemannian geometry that brings a direct visualization of the open spheres of abstract Riemannian metrics.
				
		\subsection{Preliminaries}
		We start with recalling some important notions that have been used in this article. For detailed study we recommend the book of J. M. Lee \cite{Lee2018} for Riemannian geometry and Geodesics in Riemannian manifolds, for PINN we recommend Bai \cite{Bai2022} and Liu \cite{Liu2025}.
		\subsubsection{Riemannian manifold and Riemannian metric}
		An $n$-dimensional Riemannian manifold is an $n$-dimensional smooth manifold $M$ with a globally defined 2-form $g:\chi(M)\times \chi(M)\to C^\infty(M)$ which is \begin{enumerate}
			\item positive definite i.e., $g(X,X)\ge 0,\forall X\in \chi(M)$ and $g(X,X)=0\iff X=0$.
			\item symmetric i.e., $g(X,Y)=g(Y,X),\forall X,Y\in \chi(M)$.
			\item $C^\infty(M)$-bilinear i.e., $g(fX,Y)=g(X,fY)=fg(X,Y)$, $\forall f\in C^\infty(M), \forall X,Y\in \chi(M)$.
		\end{enumerate}
		This 2-form is actually an inner product in the tangent spaces $T_p(M)$ for each $p\in M$. Let $p,q\in M$ and let $\gamma:[0,1]\to M$ be a smooth curve, then the length of the curve is given by
		\begin{equation}\label{eq_ell}
			\ell(\gamma) = \int_{0}^{1}\sqrt{g(\dot{\gamma}(t),\dot{\gamma}(t))}dt.
		\end{equation}
		This length functional is independent of parametrization. So we can measure the distance between any two points $p,q$ on the manifold $M$ by the metric induced by the length structure given by
		\begin{equation}\label{eq_metric}
			d_\ell (p,q)=\underset{\gamma}{\inf}\left\{ \ell(\gamma)|\gamma:[0,1]\overset{C^\infty}{\longrightarrow}M,\gamma(0)=p,\gamma(1)=q \right\}.
		\end{equation}
		As the distance metric $d_\ell$ depends solely on the coefficients of $g$, so $g$ is called the Riemannian metric. If one fixes an orthonormal basis $\{e_i:i=1,2,\cdots,n\}$ for $\chi(M)$ then we have a local expression (also known as local form) of $g$ given by $$g_{ij}=g(e_i,e_j).$$ These $g_{ij}$ coefficients are real valued functions of $n$ variables $x_1,x_2,\cdots, x_n$ and it forms a symmetric matrix which we denote as $g$ (the same notation) and we write
		\begin{equation}
			g = \begin{pmatrix}
				g_{11} & g_{12} & g_{13} & \cdots & g_{1n}\\
				g_{21} & g_{22} & g_{23} & \cdots & g_{2n}\\
				g_{31} & g_{32} & g_{33} & \cdots & g_{3n}\\
				\vdots & \vdots & \vdots & \ddots & \vdots\\
				g_{n1} & g_{n2} & g_{n3} & \cdots & g_{nn}
			\end{pmatrix}
		\end{equation}
		For our numerical computation we will be using this form to consider our matrix of the Riemannian metric. For any two points $p,q\in M$, a curve $\gamma:[0,1]\to M$ joining them is said to be the shortest path if $d_\ell(p,q)=\ell(\gamma)$. In the next section we will be discussing geodesics and we will show how the shortest paths becomes the geodesics.
		\subsubsection{Geodesics in Riemannian manifolds}
		
		A Riemannian metric contains important information about the geometric structure of the underlying manifold, particularly its curvature. The curvature of a Riemannian manifold is determined by its Levi-Civita connection, whose local representation is given by the Christoffel symbols. These symbols are computed from the metric tensor and its first-order partial derivatives.
		
		The \textbf{Christoffel symbols} associated with a Riemannian metric
		$g$ are defined in local coordinates by
		\begin{equation}\label{eq_christoffels}
			\Gamma^i_{jk}
			= \frac{1}{2}g^{il}
			\left(
			\frac{\partial g_{kl}}{\partial x^j}
			+ \frac{\partial g_{jl}}{\partial x^k}
			- \frac{\partial g_{jk}}{\partial x^l}
			\right),
		\end{equation}
		where $g^{il}$ denotes the components of the inverse metric tensor. The Christoffel symbols describe the local behavior of the Levi-Civita connection and are used to express covariant derivatives in local coordinates. In particular, they determine the Riemann curvature tensor, which measures the failure of covariant derivatives to commute.
		
		The Riemann curvature tensor is defined by
		\begin{equation}
			R(X,Y)Z
			= \nabla_X\nabla_Y Z
			- \nabla_Y\nabla_X Z
			- \nabla_{[X,Y]}Z,
			\qquad \forall X,Y,Z\in\chi(M),
		\end{equation}
		where $\chi(M)$ denotes the space of smooth vector fields on $M$, $\nabla$ is the Levi-Civita connection, and $[X,Y]$ is the Lie bracket of vector fields. The Riemann curvature tensor is a tensor field of type $(1,3)$, and its local components, under the convention
		$R(\partial_k,\partial_l)\partial_j
		= R^i{}_{jkl}\partial_i$, are given by
		\begin{equation}\label{eq_RiemannCT}
			R^i{}_{jkl}
			= \frac{\partial \Gamma^i_{lj}}{\partial x^k}
			- \frac{\partial \Gamma^i_{kj}}{\partial x^l}
			+ \Gamma^i_{km}\Gamma^m_{lj}
			- \Gamma^i_{lm}\Gamma^m_{kj}.
		\end{equation}
		
		Contracting the first and third tensor indices yields the Ricci curvature tensor. Its local expression is
		\begin{equation}\label{eq_RicciCT}
			R_{jk} = R^i{}_{jik},
		\end{equation}
		where the repeated index $i$ is summed over as per Einstein's summation convention. A further contraction with respect to the indices $j$ and $k$ gives the scalar curvature:
		\begin{equation}\label{eq_scalarC}
			R = g^{jk}R_{jk}.
		\end{equation}
		These local expressions are useful for computing the components of the Christoffel symbols, the Riemann curvature tensor, the Ricci tensor, and the scalar curvature symbolically.
		
		The Christoffel symbols also play a fundamental role in determining the geodesics of a Riemannian manifold. Geodesics are curves whose tangent vectors are parallel transported along the curves themselves. Let $\gamma(t)=(x^1(t),\ldots,x^n(t))$ be a smooth curve on $M$,
		with tangent vector
		\begin{equation}
			\dot{\gamma}(t)
			= \frac{dx^i}{dt}\frac{\partial}{\partial x^i}.
		\end{equation}
		A curve is a geodesic if its covariant acceleration vanishes, that is,
		\begin{equation}
			\nabla_{\dot{\gamma}}\dot{\gamma}=0.
		\end{equation}
		Expressing this condition in local coordinates and using the Christoffel symbols gives the geodesic equations
		\begin{equation}\label{eq_geodesic}
			\frac{d^2x^i}{dt^2}
			+ \Gamma^i_{jk}
			\frac{dx^j}{dt}
			\frac{dx^k}{dt}
			= 0,
			\qquad i=1,2,\cdots,n.
		\end{equation}
		Thus, the Christoffel symbols determine how the coordinate components of the velocity of a geodesic change along the curve. Since the Christoffel symbols are determined by the metric and its first-order derivatives, the geodesic equations describe how the metric governs the local behavior of geodesics. Consequently,
		both curvature computations and geodesic calculations can be carried out symbolically once the metric tensor is specified.
		
		\subsubsection{Physics-Informed Neural Network}
		Physics-Informed Neural Network is a particular type of neural network from Deep Learning which helps us to solve mathematical problems by incorporating known physical laws, such as differential equations, into the training process. In ordinary machine learning, a neural network generally learns from a given set of data. In a PINN, the network learns from data and from the governing equations of the problem. Some times we often let the system learn from the governing equations only to produce a data independent model and then compare the model's solution to known experimental data for accuracy. This makes PINNs particularly useful for solving ordinary differential equations (ODEs), partial differential equations (PDEs), and certain inverse problems. The framework was popularized by Raissi, Perdikaris, and Karniadakis in their 2019 paper \cite{Raissi2019}.
		\\
		\vskip6pt
		\noindent\textbf{Mathematical formulation of PINN:}\\
		Consider a general differential equation of the form
		\begin{equation}\label{eq_pde}
			\mathcal{D}[u](x)=f(x), \qquad x\in\Omega,
		\end{equation}
		where $\Omega\subset\mathbb{R}^n$ denotes the computational domain, $u$ is the unknown function, $\mathcal{D}$ is a differential operator, and $f$ is a given source term.
		
		A Physics-Informed Neural Network (PINN) approximates the unknown solution $u$ by a neural network $u_\theta$, where $\theta$ represents the trainable parameters of the network. Thus, we write
		\begin{equation}\label{eq_nn_approx}
			u(x)\approx u_\theta(x).
		\end{equation}
		
		To incorporate the governing differential equation into the training process, we define the \textit{physics residual} as 
		\begin{equation}\label{eq_physics_residual}
			p_\theta(x)
			=\mathcal{D}[u_\theta](x)-f(x).
		\end{equation}
		If the neural network satisfies the differential equation exactly, then $p_\theta(x)=0$ throughout the domain $\Omega$. In practice, the residual is evaluated at a finite collection of collocation points $\{x_i\}_{i=1}^{N_p}\subset\Omega$, where $N_p$ denotes the number of collocation points. The corresponding physics loss is defined by
		\begin{equation}\label{eq_physics_loss}
			\mathcal{L}_{p}
			=\frac{1}{N_p}\sum_{i=1}^{N_p}
			\left|p_\theta(x_i)\right|^2.
		\end{equation}
		
		Suppose that the differential equation is subject to the boundary condition
			\begin{equation}\label{eq_boundary_condition}
			u(x)=\phi(x), \qquad x\in\partial\Omega,
		\end{equation}
		where $\partial\Omega$ denotes the boundary of the domain and $\phi$ is a prescribed boundary function. The boundary condition is incorporated into the training process through the boundary loss
		\begin{equation}\label{eq_boundary_loss}
			\mathcal{L}_{b}
			=\frac{1}{N_b}\sum_{i=1}^{N_b}
			\left|u_\theta(x_i^b)-\phi(x_i^b)\right|^2,
		\end{equation}
		where $\{x_i^b\}_{i=1}^{N_b}$ are points on the boundary $\partial\Omega$, and $N_b$ is the number of boundary points.
		
		The total loss function is then constructed as a weighted combination of the physics loss and the boundary loss:
		\begin{equation}\label{eq_total_loss}
			\mathcal{L}(\theta)
			=\lambda_p\mathcal{L}_{p}
			+\lambda_b\mathcal{L}_{b},
		\end{equation}
		where $\lambda_p$ and $\lambda_b$ are nonnegative weighting parameters that control the relative contributions of the two loss terms. The optimal network parameters are obtained by solving the minimization problem
		\begin{equation}\label{eq_optimization}
			\theta^*
			=\operatorname*{arg\,min}_{\theta}\mathcal{L}(\theta).
		\end{equation}
		
		The minimization is typically performed using gradient-based optimization algorithms. The derivatives of the neural network output with respect to the input coordinates, which are required to evaluate the differential operator $\mathcal{D}$, can be computed efficiently using automatic differentiation.
		
		If observational data are available, an additional data-fitting term can be incorporated into the loss function. For example, given observations $\{(x_i^d,u_i^d)\}_{i=1}^{N_d}$, the data loss can be defined as
		\begin{equation}\label{eq_data_loss}
			\mathcal{L}_{d}
			=\frac{1}{N_d}\sum_{i=1}^{N_d}
			\left|u_\theta(x_i^d)-u_i^d\right|^2.
		\end{equation}
		In this case, the total loss may be written as \begin{equation}\label{eq_total_loss_data}
			\mathcal{L}(\theta)
			=\lambda_p\mathcal{L}_{p}
			+\lambda_b\mathcal{L}_{b}
			+\lambda_d\mathcal{L}_{d},
		\end{equation}
		where $\lambda_d$ is the weight assigned to the data loss.
		
		Thus, a PINN combines the governing differential equation, boundary conditions, and, when available, observational data within a unified optimization framework. This approach allows the neural network to approximate solutions of differential equations even when exact solutions or extensive training data
		are unavailable.
		
		In our case, we are going to implement PINN into solving geodesic equations \eqref{eq_geodesic} formed by the metric $g$ in \eqref{eq_metric} with suitably chosen initial conditions. It should be mentioned here that we are not going to use experimental data for our PINN model, so that the network will learn the system by itself and we will use the \texttt{SciPy} solution only to compare the accuracy of the model's solution. We specify all the working methodologies in the next section.

		\section{Working Method and Configurations}
		\subsection{Working method}
		
		\noindent\textbf{Step-I (Dimension):} First we fix the dimension $n$. For our case it is $n=3$. We are taking a 3 dimensional manifold so that we can have a relatable visualization. Jumping into higher dimension will be easy once we can see what is happening in the lower dimension. One may also consider the same experiment on dimension 2. We also intent to investigate dimension 4 with space-time metrics with possibly singularities in near future.\\
		
		\noindent\textbf{Step-II (Domain and Coordinate):} Next we fix a domain $\Omega\subset \mathbb{R}^3$ with coordinate $x,y,z$ in which the Riemannian metric will be defined. In general case the unit cube $[0,1]^3$ serves the purpose. But in our case we relax this case by choosing globally defined metric which has no singularity in the whole of $\mathbb{R}^3$, thus in our case we specify $\Omega=\{(x,y,z)\in\mathbb{R}^3:x>0\}\subset \mathbb{R}^3$. However, the parameter $\tau$ that controls the geodesic trajectory is taken in the interval $[0,2]$ so that half length of the geodesic can be traced out in $[0,1]$. The reason behind choosing $x>0$ is specified in the next section.\\
		
		\noindent\textbf{Step-III (The Riemannian metric):} After domain specification we specify the Riemannian metric $g$ in matrix form. In our case we are dealing with a 3 dimensional matrix which is given by
		$$g=\begin{pmatrix}
			1 & 0 & 0\\
			0 & x^2 & 0\\
			0 & 0 & 1
		\end{pmatrix}, x>0.$$
		We prefer to start with a simple metric and just changing one element may seem to be not so interesting but the results are very impressive as this will directly show us the effect of changing one parameter in the metric or in the line element. In our case the line element is $$ds^2=dx^2+x^2dy^2+dz^2, x>0.$$	We are going to visualize the open spheres under this metric $g$. We keep the metric to be positive definite for now but in future we may extend this to negative definite or indefinite metrics as well.\\
		
		\noindent\textbf{Step-IV (Study the geometry of $g$):} Now we use symbolic calculations to find the expressions for Christoffel symbols using \eqref{eq_christoffels}, Riemann curvature tensor using \eqref{eq_RiemannCT}, Ricci curvature using \eqref{eq_RicciCT}, Scalar curvature using \eqref{eq_scalarC}.\\
		
		\noindent\textbf{Step-V (Einstein manifold and Flatness checking) [Optional]:}	It is not a necessary step but we like to keep it in the algorithm part so that we can report this additional information. A manifold is said to be Einstein if $Ric=\lambda g$ for some scalar $\lambda$. A manifold is said to be flat if its Riemann curvature vanishes identically, i.e., all of the Riemann curvature components are zero. Similarly, a Ricci flat manifold is the one with vanishing Ricci curvature. So Einstein manifolds can be considered as an extension of Ricci flat manifolds. It has various applications in Astrophysics and Cosmology. We can check for all nonzero terms of Riemann and Ricci tensors after the computation. Depending on their vanishing we will present our verdict. For Einstein manifold checking, we simply compute $Ric-\lambda g$ and compare it with the $O_{3\times 3}$ (the null matrix), here we solve for $\lambda$ and depending on the number of solutions we will say whether it is Einstein manifold or not.\\
		
		\noindent\textbf{Step-VI (Geodesic equations):} Next using the Christoffel symbols and the coordinates, we form the geodesic equations. We solve them using \texttt{SciPy}'s integrate module. In our case we will have the solution of the form $G = \{(x(\tau),y(\tau),z(\tau)):\tau \in [0,2]\}$. This will form the geodesics in a parametric form. Here $x(\tau),y(\tau),z(\tau)$ satisfies the geodesic equations \eqref{eq_geodesic} for our specified metric $g$.\\
		
		\noindent\textbf{Step-VII (PINN):} Next we implement PINN to solve the equations and compare the results with \texttt{SciPy}'s outcome. After running the optimizer for 10 different seeds we calcualte the overall root mean square error (RMSE) using this formula
		$$
		\text{Overall RMSE}
		=
		\sqrt{
			\frac{1}{3N}
			\sum_{i=1}^{N}
			\sum_{j=1}^{3}
			\left(
			\hat{u}_{ij}-u_{ij}
			\right)^2
		}
		$$
		
		where
		\begin{eqnarray*}
			N &=& 1000, \quad \text{number of sampled points},\\
			\hat{u}_{ij} &=& \text{PINN-predicted value of coordinate } j \text{ at point } i,\\
			u_{ij} &=& \text{numerical reference value of coordinate } j \text{ at point } i,\\
			i &=& 1,2,\cdots,N,\\
			j &=& 1,2,3 \quad \text{corresponding to } x,y,z.\\
		\end{eqnarray*}
		Please refer to Table \ref{tab:pinn_error_metrics} for the numerical values.\\
		
		\noindent\textbf{Step-VIII (Spike Diagram):} After that we implement our custom made python code to plot the geodesic spike diagram. The spike diagram is basically tracing all the points on the surface of a sphere (finite number of points) from the center of a sphere derived from the geodesic equations. In the later section we will discuss in detail about this diagram.\\
		
		\noindent\textbf{Step-IX (Mesh creation for the geodesic open sphere):} Using python's inbuilt codes we are able to construct the surface formed by the end points of the geodesics derived in the earlier step. We call this closed surface a geodesic sphere or geodesic open sphere.\\
		
		\noindent\textbf{Step-X (Analysis on the newly formed geodesic sphere):} Now we do certain analysis by answering questions like: \textit{How much it is different from the Euclidean sphere? Which side of the geodesic sphere is compressed and which side is expanded compared to the Euclidean sphere? What is the relative volume / surface area difference between them?}

		\subsection{Configurations}\label{sec_dev_config}
		\subsubsection{Device Configurations}
		We are using free version of Google Colab (Python 3 Google Compute Engine backend (GPU)) T4 GPU with system RAM 12.7 GB, GPU RAM 15 GB, Disk 112.6 GB. However, the code can be executed through Python 3 CPU, but it would take a little long time for the whole evaluation.
		\subsubsection{Architecture of the Physics-Informed Neural Network}
		Our Physics-Informed Neural Network consists of 6 layers in total in which 1 input layer (with 1 neuron $\tau$), 4 hidden layers each equipped with 64 neurons each and 1 output layer (with 3 neurons $x(\tau),y(\tau),z(\tau)$). Only two types of derivatives are involved $\frac{d}{d\tau}$ and $\frac{d^2}{d\tau^2}$ and they are used under automatic differentiation. Three ODEs are involved with 6 initial conditions (3 conditions for initial points and 3 conditions for initial velocities foe each of the output variables). We used \textbf{Adam} optimizer with \textbf{learning rate} $10^{-3}$ and number of \textbf{iterations} are $10000$ for each execution. Our loss function is given by $$\mathcal{L}(\theta) = \lambda_p \mathcal{L}_p+\lambda_b \mathcal{L}_b,$$ where $\mathcal{L}_p$ denotes the loss in physics (i.e., loss in learning the governing equations) and $\mathcal{L}_b$ is the loss in learning the boundary/initial conditions. The $\lambda_p,\lambda_b$ are the weights for the physics loss and boundary loss functions. We are letting the network learn the ODEs by the physics and boundary/initial conditions only, and we do not invoke any experimental data points. We only use the \texttt{SciPy} data to compare the results with the PINN. We are executing the neural network for $10$ random seeds starting from $1234$ and up to $1243$. This makes the network add different weights to the loss function and from here we find the best seed by choosing the minimum overall root mean square error (RMSE).
		We used \textbf{Tanh} as activation function and \textbf{Glorot Uniform} as initializer.
		The network architecture diagram is given here for readers convenience.
		\begin{figure}[H]
			\centering
			\caption{PINN Architecture under Adam optimizer}
			\label{fig:pinnarchitecture}
			\includegraphics[width=0.9\linewidth]{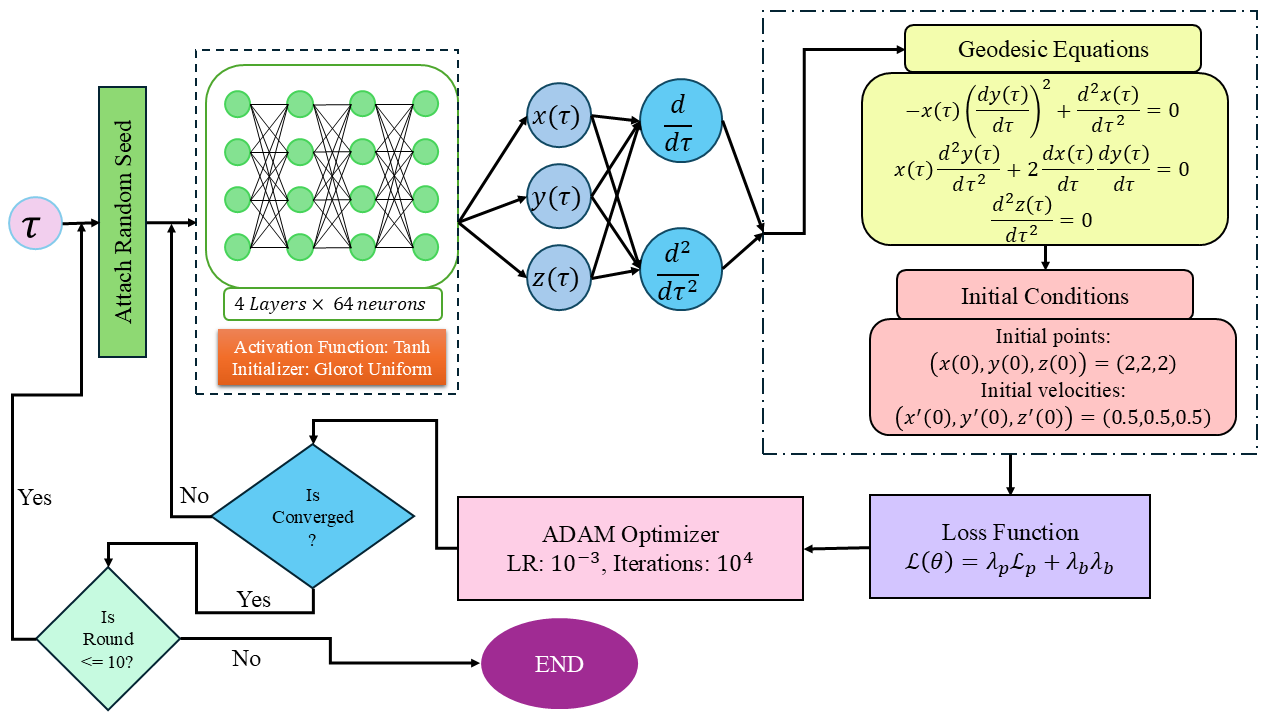}
			
		\end{figure}
		
		\section{Results and Discussions}
		\subsection{Computation time}
		First let us mention the computation time for the execution of the Python script on Google Colab T4 GPU using the same configuration as mentioned in Section \ref{sec_dev_config}.
		\begin{itemize}
			\item Calculating symbolic Christoffel symbol, Riemann, Ricci and Scalar curvature, Einstein manifold checking and Flatness checking: \textbf{6.220 sec}
			\item Finding geodesic equations: \textbf{0.240 sec}
			\item Installation of DeepXDE: \textbf{10.100 sec}
			\item PINN 10 seeds each with 10k iterations total time: \textbf{315.777 sec}
			\item Volume comparison : \textbf{60.05 sec}
		\end{itemize}
		Estimated Total time: \textbf{6 minutes 32.387 seconds.}
		\begin{remark}
			It should be mentioned here that the parameters we have chosen for this experiment are not random. They are chosen after multiple execution of the program and several failures. Sometimes the method took long time to respond and sometimes it was faster but the errors were drastic. After a careful observation and keeping the time of execution in mind we have finalized the parameter values.
		\end{remark}
		\subsection{Results on the geometry of the metric}
		We are considering $(M,g)$ as our manifold where $M=\{(x,y,z)\in \mathbb{R}^3, x>0\}$, $x,y,z$ as the coordinates in the manifold and the metric $g$ is defined by 
		$$g=\begin{pmatrix}
			1 & 0 & 0\\
			0 & x^2 & 0\\
			0 & 0 & 1
		\end{pmatrix}.$$ So the line element becomes
		$$ds^2 = dx^2+x^2dy^2+dz^2.$$
		\begin{remark}
			We have chosen this region so that $g$ remains positive definite, which is necessary for our experiment. Also the reason behind choosing such simple metric is that we had no idea about the outcome. To make sure we have very good control over the metric we are changing just one diagonal element and keeping the remaining the same as Euclidean metric. In future, we have plan to investigate over complex metrics, degenerate and indefinite metrics along with metrics with singularity. 
		\end{remark}
				
		\noindent The nonzero components of the Christoffel symbols for $g$ are
		\begin{eqnarray*}
			\Gamma^{1}_{22} = -x\\
			\Gamma^{2}_{12} = \Gamma^{2}_{21}= \frac{1}{x}.\\
		\end{eqnarray*}
		All of the Riemann, Ricci and Scalar curvature for the metric $g$ are zero. So the manifold is a \textbf{flat manifold}.\\
		The geodesic equations for $g$ are given by
		\begin{equation}\label{eq_geodesic_g}
			\begin{cases}
				-x(\tau)\left(\frac{d}{d\tau}y(\tau)\right)^2+\frac{d^2}{d\tau^2}x(\tau)=0,\\
				x(\tau)\frac{d^2}{d\tau^2}y(\tau)+2\frac{d}{d\tau}x(\tau)\frac{d}{d\tau}y(\tau)=0,\\
				\frac{d^2}{d\tau^2}z(\tau)=0.
			\end{cases}
		\end{equation}
		To solve the system \eqref{eq_geodesic_g} we assume that at $\tau\in[0,2]$ and at $\tau=0$ we have
		\begin{eqnarray*}
			\begin{cases}
				x(0)=2,\\
				y(0)=2,\\
				z(0)=2,\\
			\end{cases} &,&
			\begin{cases}
				x'(0)=1/2,\\
				y'(0)=1/2,\\
				z'(0)=1/2.\\
			\end{cases}
		\end{eqnarray*}
		
		\subsection{Results of PINN}
		Now we apply both \texttt{SciPy} Integrate Solve IVP with \textbf{DOP853} method (which is 8th order Runge-Kutta Numerical Solver with 5th order error estimator) and our PINN to get the solution. We have done random seed experiments starting seen number from 1234 up to 1243 and the following graphs shows the test-train loss for each experiment.
		
		\begin{figure}[H]
			\centering
			\caption{Test loss and Training loss graph for 10 random seed PINN execution}
			\label{fig:testloss}
			\includegraphics[width=1\linewidth]{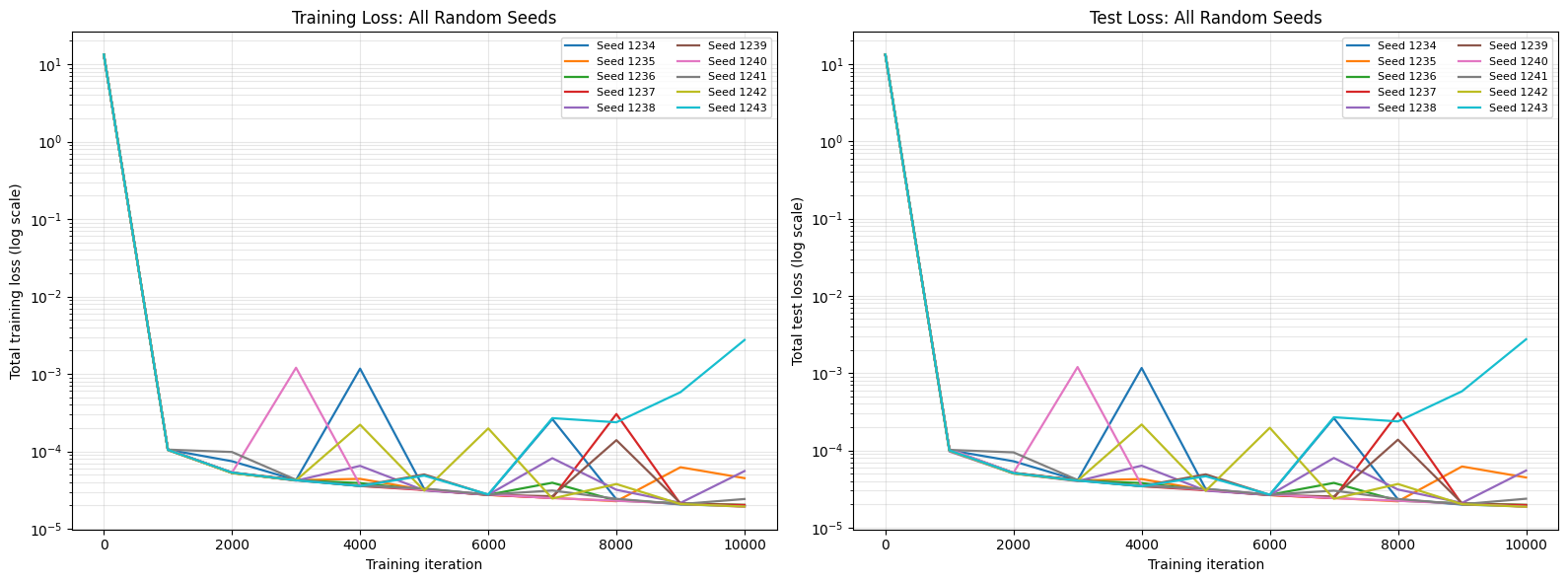}
		\end{figure}
		The Table \ref{tab:pinn_error_metrics} shows the error metrics of each experiment. As you can see the best performance on the basis of minimum overall RMSE is achieved by the first seed \textbf{1234} (See Table \ref{tab:pinn_error_metrics}). So interested reader can use this seed to reconstruct the same result or may find another seed for better outcome. The Table \ref{tab:final_position} gives a comparison of final position at $\tau=2$ predicted by Numerical solution and PINN solution. Tables are placed at the end of the section for better readability.
		
\begin{figure}[H]
			\centering
			\caption{}
			\label{fig:pinnxyzvstau}
			\includegraphics[width=0.7\linewidth]{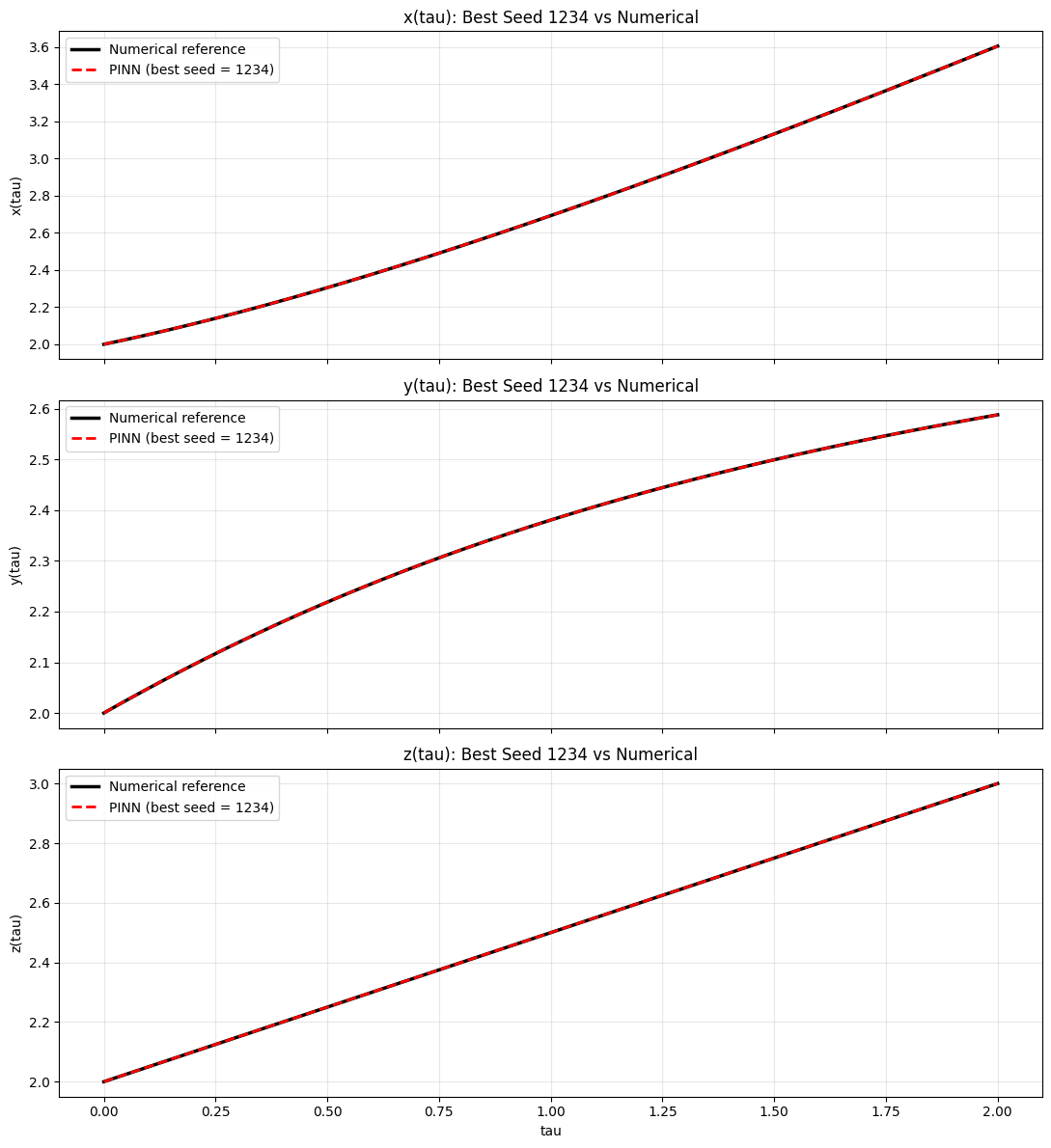}
		\end{figure}

		Now let us see the solution comparison between PINN and Numerical method. See Figure \ref{fig:pinnxyzvstau}. The PINN solution lines almost coincides with te numerical solution.
		
		To understand how close is the PINN prediction we show the box plots of the errors generated at each seed.
		
		\begin{figure}[H]
			\centering
			\caption{}
			\label{fig:errorboxplotallruns}
			\includegraphics[width=\linewidth]{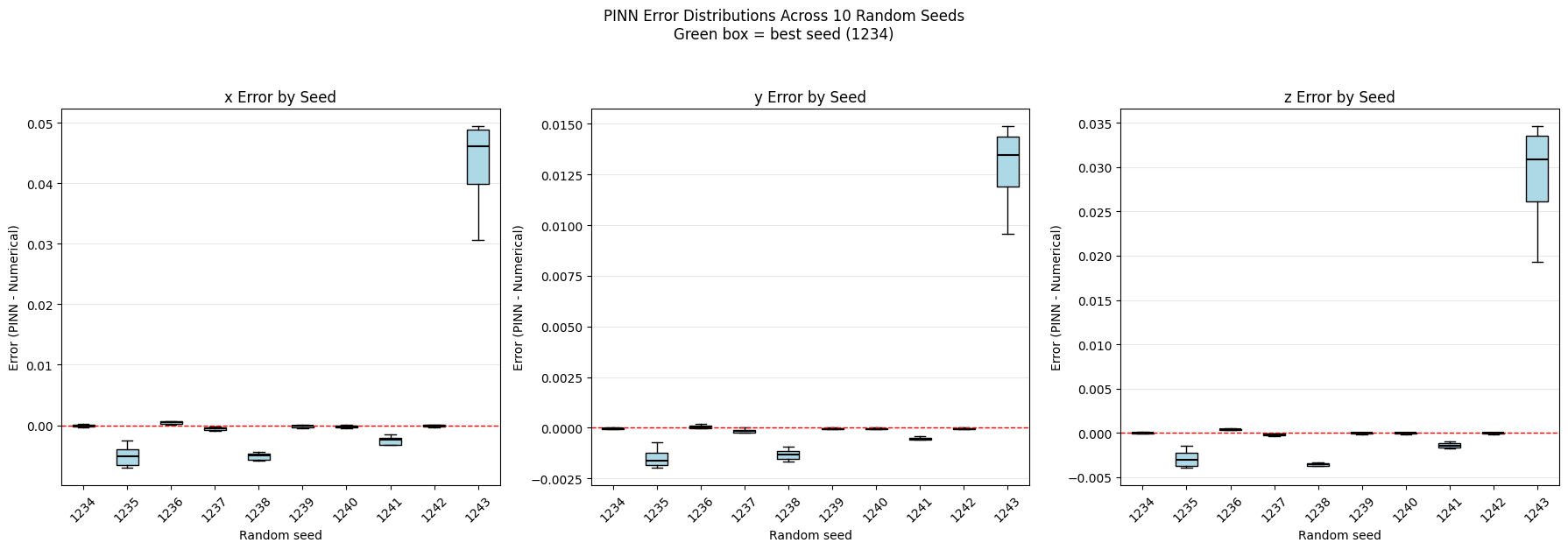}
		\end{figure}
		Honestly, we did not expected the model to produce such a small box plot for each run. It is very good to get such small errors for each of the run except for the last one. This shows the model's stability over multiple runs. We suspect that further tuning the model's configuration may lead to a more stabilized prediction. For this run the overall RMSE is $1.144660\times 10^{-4}$. If one runs the code again, they may get a different RMSE and different seed.
	
		\subsection{Results on Geodesic Spike Diagram} Now we are ready to visualize the geodesic spike diagram. Let $p$ be a point fixed in $\mathbb{R}^3$ and let $S(p,1)$ be the unit sphere around the point $p$ i.e., $S(p,1)=\{x\in \mathbb{R}^3:||x-p||_{2}=1\}$. First we discretize the sphere into finite number of points. Then draw the geodesics (produced by $g$) starting from the center $p$ to all the other points on the surface. This will create a spike like diagram. For the Euclidean metric, the shortest paths are straight lines. If we draw the same for Euclidean metric then all we are going to see a bunch of finite sticks glues at one point. Topologically speaking it is equivalent to $[0,1]^{100}/\sim$ where $\sim$ identifies the $0$ of each interval. This is what we mean by spike diagram. For arbitrary Riemannian metric we are interested is the diagram, this will give us an idea what does the unit sphere will look like under the metric $g$.
		
		For this visualization we are creating $15$ points along the $\theta$ coordinate and $30$ points in the $\phi$ coordinate on the unit sphere. As a total we have $450$ points in the surface to which geodesics are going to be connected. 
		
		However, \textbf{we will see $422$ geodesics in the plot}. Why? To visualize the geodesics emanating from the center (fixed point), we discretize the unit sphere of initial directions using $N_\theta=15$ polar subdivisions and $N_\phi=30$ azimuthal subdivisions. Since all azimuthal angles represent the same direction at each pole, only one direction is retained at each pole. Hence, the total number of distinct initial directions is 
		$$
		N_{\mathrm{geo}}=(N_\theta-1)N_\phi+2=(15-1)\times 30+2=422.
		$$
		Each initial direction determines a geodesic initial-value problem, which is numerically solved using \texttt{SciPy}'s DOP853 method. Thus, at most 422 geodesic trajectories are obtained, depending on the success of numerical integration. This number depends on the chosen angular resolution and does not represent the total number of geodesics on the manifold.
		
		\begin{figure}[H]
			\centering
			\caption{Geodesic Spike diagram for the metric \eqref{eq_metric}}
			\label{fig:geodesicspikediagram}
			\includegraphics[width=0.7\linewidth]{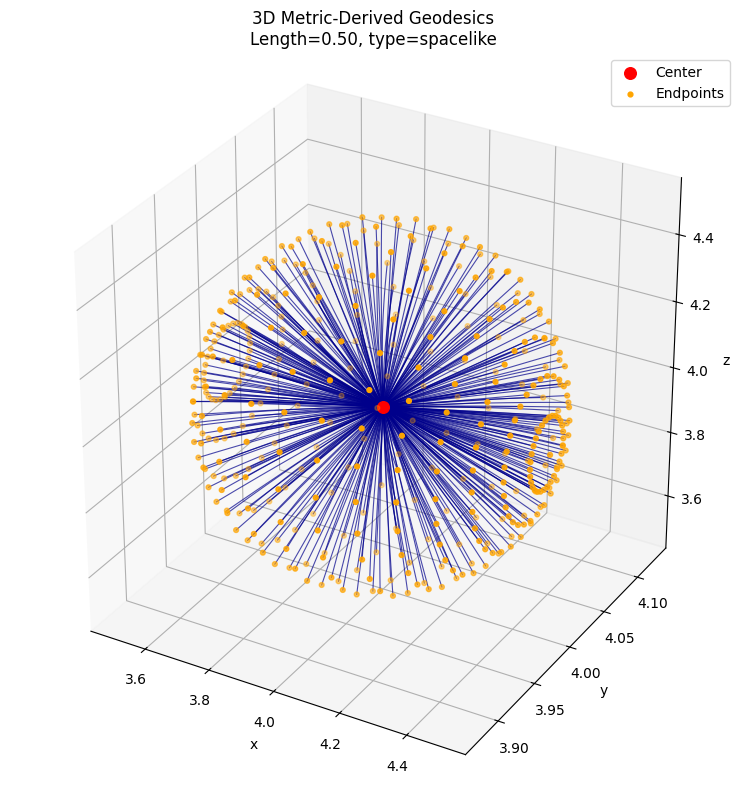}
		\end{figure}
		
		\subsection{Construction of Geodesic Sphere}
		Now we make a triangular mesh to construct the surface formed by the displaced points of the deformed sphere. We call this the \textbf{geodesic sphere}. Since the metric we have considered in \eqref{eq_metric} is very close to Euclidean metric, so plotting the open sphere is of no use as the human eye can not tell the difference between this image and the Euclidean sphere. To explicitly visualize the change we compare the newly generated geodesic sphere with Euclidean unit sphere $\mathbb{S}^2$. First comparison comes in form of compression and decompression. Let us define the displacement at the $i$-th point as
		
		$$\Delta_i =\left\|\mathbf{x}_i^{\text{geodesic}}-\mathbf{x}_i^{\text{Euclidean}}\right\|_2$$
		
		where $\mathbf{x}_i^{\text{geodesic}}$ is the endpoint of the $i$-th geodesic, $\mathbf{x}_i^{\text{Euclidean}}$ is the corresponding point on the reference Euclidean sphere, $\Delta_i$ is the Euclidean coordinate distance between the two corresponding points. Plotting the heat map on the geodesic sphere gives us a nice visualization. 
		
		\begin{figure}[H]
			\centering
			\caption{Geodesic sphere created from the geodesic spike diagram}
			\label{fig:heatmapgeodesicsphere}
			\includegraphics[width=0.7\linewidth]{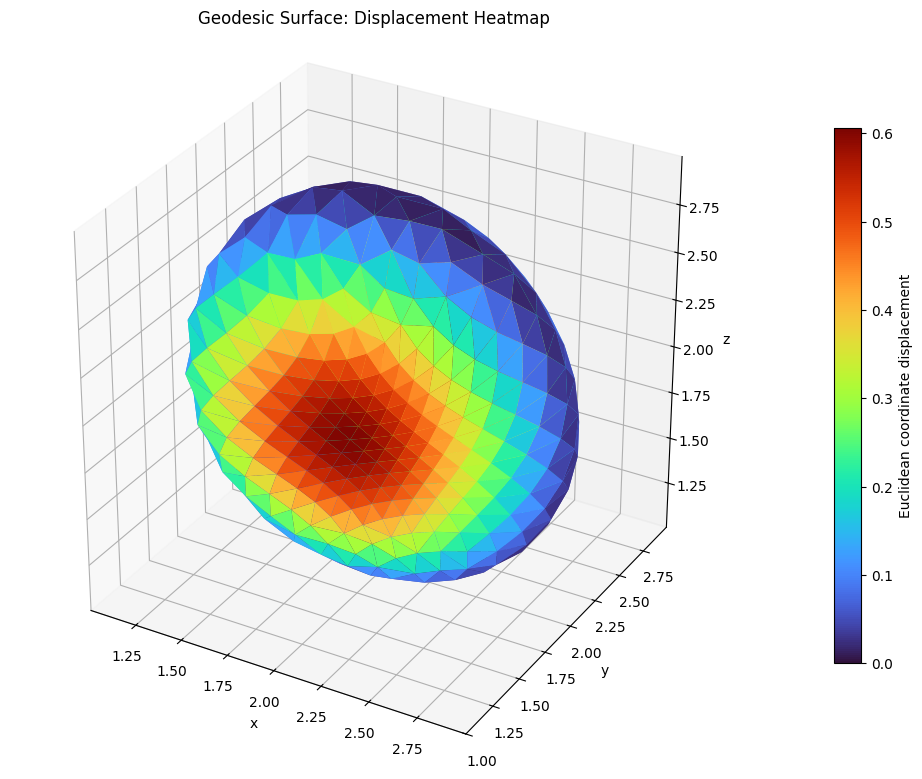}
		\end{figure}
		
		\subsection{Spatial distance comparison}
		This graph (Figure \ref{fig:heatmapgeodesicsphere}) gives an idea about the parts of the open sphere of the Riemannian metric where it has a difference with the Euclidean unit sphere. Next we consider the signed radial deviation. The signed radial deviation measures how far a geodesic endpoint lies from the reference sphere in the radial direction. It is defined as
		
		$$\delta_i =\left\|\mathbf{x}_i^{\text{geodesic}}-p\right\|_2 - R$$
		
		where $p=(2,2,2)$ is the common center of the two constructions and $R$ is the radius of the reference Euclidean sphere, here $R=1$. The sign of $\delta_i$ indicates the endpoint's radial position:
		\begin{enumerate}
			\item $\delta_i > 0$: the geodesic endpoint lies outside the reference sphere.
			\item $\delta_i < 0$: the geodesic endpoint lies inside the reference sphere.
			\item $\delta_i = 0$: the geodesic endpoint lies on the reference sphere.
		\end{enumerate}
		\begin{remark}
			The displacement $\Delta_i$ measures the full Euclidean distance between corresponding points, whereas $\delta_i$ measures only the difference in their distances from the common center. These quantities are generally different. Neither directly represents intrinsic curvature.
		\end{remark}
		
		Before going into the signed distance deviation let us see the effect of $\Delta_i$ in the Euclidean sphere.
		
		\begin{figure}[H]
			\centering
			\caption{Displacement heat maps ($\Delta_i$) on $\mathbb{S}^2$}
			\label{fig:distancedeviationsphere}
			\includegraphics[width=0.7\linewidth]{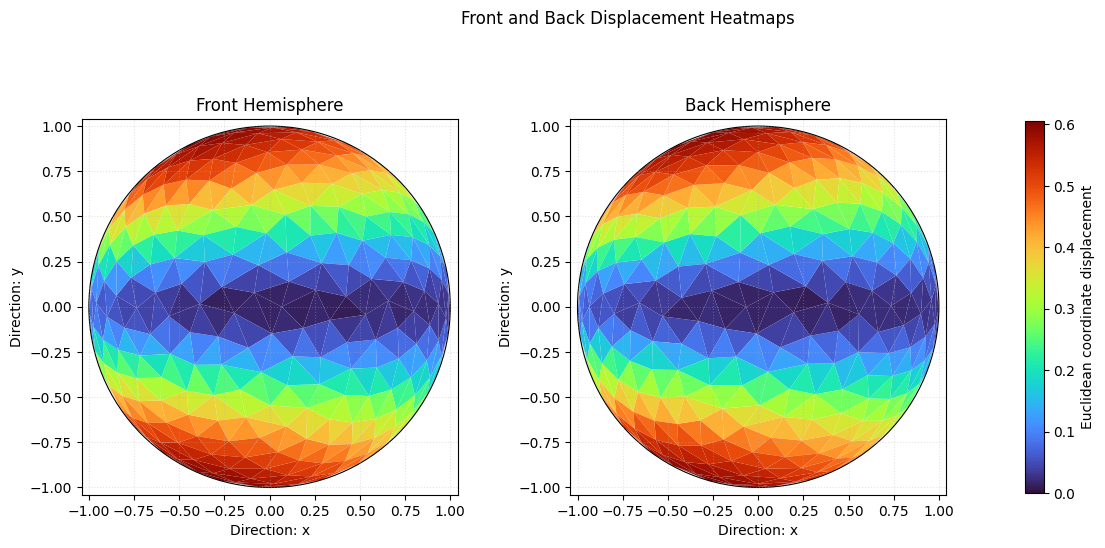}
		\end{figure}
		
		This image is on the Euclidean sphere, and this represents which part has a change when the metric is invoked. Now for the signed distance deviation heatmap on the Euclidean sphere we present the following diagram.
		
		\begin{figure}[H]
			\centering
			\caption{Signed distance deviation heat maps ($\delta_i$) on $\mathbb{S}^2$}
			\label{fig:signeddistancedeviationsphere}
			\includegraphics[width=0.7\linewidth]{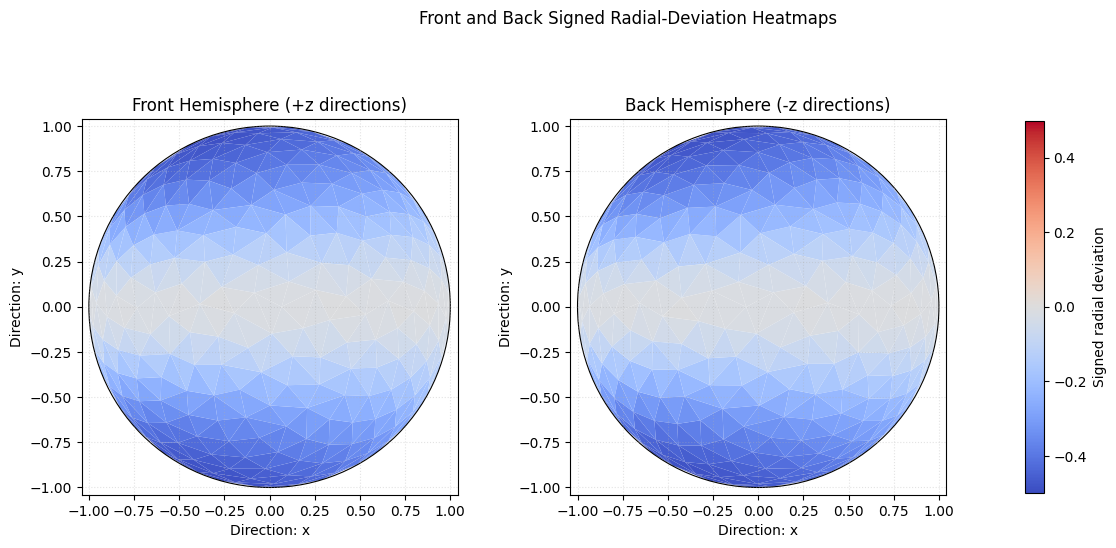}
		\end{figure}
		
		It can be seen from here that contraction is happening when the metric is involved in $\mathbb{R}^3$. This is what we are looking for, this visualization answers one of the fundamental questions which is ``\textbf{what does the Riemannian metric does to the unit sphere?}" Answer: Figure \ref{fig:signeddistancedeviationsphere} is what is does to the sphere.
		
		The next question comes to mind is, \textbf{Can we compare the change in volume and surface for this newly generated geodesic sphere with the Euclidean sphere?} Yes we can and now we are going to do the same.
		
		Since the surface of the geodesic sphere is made of triangular meshes so we can numerically compute all the necessary information.
		
		\subsection{Intrinsic Riemannian Area and Volume Comparison}
		
		Till now we have $(M,g)$ as a three-dimensional Riemannian manifold with metric tensor $g=(g_{ij})$. For a geodesic sphere of radius $R$ centered at $p\in M$, we compare its intrinsic area and volume with those of a Euclidean sphere and sphere of the same radius. Unlike Euclidean mesh calculations, which use the vertex
		coordinates and triangular faces, intrinsic Riemannian quantities are determined by the metric tensor $g$.
		
		The intrinsic area of the geodesic sphere is calculated using the induced metric:
		$$
			A_g(R)=\int_0^\pi\int_0^{2\pi}
			\sqrt{EH-F^2}\,d\phi\,d\theta,
		$$
		where
		$$
		E=g(X_\theta,X_\theta),\qquad
		F=g(X_\theta,X_\phi),\qquad
		H=g(X_\phi,X_\phi),
		$$
		and $X(\theta,\phi)$ parameterizes the geodesic sphere.
		
		The intrinsic volume of the geodesic sphere is calculated by integrating the Riemannian volume density:
		$$
			V_g(R)=\int_0^R\int_0^\pi\int_0^{2\pi}
			\sqrt{\det g(F)}\,|\det J_F|\,
			d\phi\,d\theta\,ds,
		$$
		where $F(s,\theta,\phi)$ is the geodesic polar map and
		$$
		J_F=
		\frac{\partial(x^1,x^2,x^3)}
		{\partial(s,\theta,\phi)}
		$$
		is its coordinate Jacobian matrix. These formulas assume that the geodesic polar map is regular and injective over the sphere.
		
		For comparison, the Euclidean reference values are $(R=1)$,
		$$
		A_{\mathrm{Euc}}=4\pi R^2,
		\qquad
		V_{\mathrm{Euc}}=\frac{4}{3}\pi R^3.
		$$
		
		The percentage differences are
		$$
		\Delta_A=
		100\times\frac{A_g(R)-A_{\mathrm{Euc}}}
		{A_{\mathrm{Euc}}}\%,
		$$
		$$
		\Delta_V=
		100\times \frac{V_g(R)-V_{\mathrm{Euc}}}
		{V_{\mathrm{Euc}}}\%.
		$$
		
		Thus, Euclidean mesh area and volume describe geometry in ambient Euclidean coordinates, whereas intrinsic Riemannian area and volume incorporate the metric tensor $g$.
		
		\begin{table}[H]
			\centering
			\caption{Intrinsic Riemannian Geodesic sphere Comparison}
			\begin{tabular}{lr}
				\hline
				\textbf{Parameter} & \textbf{Value} \\
				\hline
				Center & $(2,2,2)$ \\
				Geodesic radius & $1.0$ \\
				Base grid $(N_\theta,N_\phi,N_r)$ & $(30,60,40)$ \\
				\hline
				\multicolumn{2}{c}{\textbf{Surface Area}} \\
				\hline
				Euclidean reference area & $12.5663706144$ \\
				Riemannian sphere area & $12.5208785515$ \\
				Difference & $-0.0454920628$ \\
				Percentage difference & $-0.362014\%$ \\
				\hline
				\multicolumn{2}{c}{\textbf{Volume}} \\
				\hline
				Euclidean reference volume & $4.1887902048$ \\
				Riemannian sphere volume & $4.1735465537$ \\
				Difference & $-0.0152436511$ \\
				Percentage difference & $-0.363915\%$ \\
				\hline
				\multicolumn{2}{c}{\textbf{Grid-Refinement Diagnostics}} \\
				\hline
				Refined grid $(N_\theta,N_\phi,N_r)$ & $(45,90,60)$ \\
				Refined sphere area & $12.5461407142$ \\
				Refined sphere volume & $4.1820115974$ \\
				Area change relative to base & $+0.201354\%$ \\
				Volume change relative to base & $+0.202416\%$ \\
				\hline
			\end{tabular}
		\end{table}
		
		This shows that at the base resolution, the computed Riemannian sphere area and sphere volume are approximately $0.362014\%$ and $0.363915\%$ smaller than their respective Euclidean reference values. Grid refinement changes the area and volume by approximately $0.201354\%$ and $0.202416\%$, respectively. These changes indicate numerical sensitivity. Interested readers may investigate on further refinement to establish the convergence. This numerical computations justifies why there are so many blue regions in the Figure \ref{fig:signeddistancedeviationsphere}.
		
		\paragraph{\textbf{Mathematical Assumptions we have taken}.}
		\begin{itemize}
			\item The metric must remain smooth and positive-definite throughout the sampled region.
			\item The exponential map must remain injective over the geodesic sphere.
			\item The coordinate Jacobian is approximated using finite differences, so the results should be checked at successively finer resolutions.
		\end{itemize}
		
		\paragraph{\textbf{Numerical limitations}}
		The computed area and volume depend on the chosen grid resolution and the numerical approximation of the coordinate Jacobian. Although grid refinement produces relatively small changes in the computed quantities, these results alone do not establish convergence. Further refinement and error analysis are required to assess the numerical accuracy of the intrinsic area and volume estimates.
		
		\begin{table}[h]
			\centering
			\caption{Error metrics for ten independent PINN runs with different random seeds.}
			\label{tab:pinn_error_metrics}
			\resizebox{\textwidth}{!}{%
				\begin{tabular}{ccrrrr}
					\toprule
					\textbf{Seed} & \textbf{Coordinate} & \textbf{MAE} & \textbf{RMSE} &
					\textbf{MaxAE} & \textbf{Overall RMSE} \\
					\midrule
					
					\multirow{3}{*}{1234}
					& $x$ & $1.475333\times10^{-4}$ & $1.801110\times10^{-4}$ & $3.399890\times10^{-4}$ & \multirow{3}{*}{$1.144660\times10^{-4}$} \\
					& $y$ & $3.439841\times10^{-5}$ & $4.341200\times10^{-5}$ & $7.044621\times10^{-5}$ & \\
					& $z$ & $6.393794\times10^{-5}$ & $7.058888\times10^{-5}$ & $1.122085\times10^{-4}$ & \\
					\midrule
					
					\multirow{3}{*}{1235}
					& $x$ & $5.113784\times10^{-3}$ & $5.295075\times10^{-3}$ & $7.063421\times10^{-3}$ & \multirow{3}{*}{$3.638967\times10^{-3}$} \\
					& $y$ & $1.526544\times10^{-3}$ & $1.569376\times10^{-3}$ & $1.985233\times10^{-3}$ & \\
					& $z$ & $2.936206\times10^{-3}$ & $3.037348\times10^{-3}$ & $3.953934\times10^{-3}$ & \\
					\midrule
					
					\multirow{3}{*}{1236}
					& $x$ & $4.181847\times10^{-4}$ & $4.715432\times10^{-4}$ & $7.168459\times10^{-4}$ & \multirow{3}{*}{$3.653336\times10^{-4}$} \\
					& $y$ & $7.321597\times10^{-5}$ & $9.054659\times10^{-5}$ & $2.007484\times10^{-4}$ & \\
					& $z$ & $4.063160\times10^{-4}$ & $4.121339\times10^{-4}$ & $4.953508\times10^{-4}$ & \\
					\midrule
					
					\multirow{3}{*}{1237}
					& $x$ & $5.677972\times10^{-4}$ & $6.236564\times10^{-4}$ & $9.560724\times10^{-4}$ & \multirow{3}{*}{$3.917300\times10^{-4}$} \\
					& $y$ & $1.602889\times10^{-4}$ & $1.783923\times10^{-4}$ & $2.461938\times10^{-4}$ & \\
					& $z$ & $1.728306\times10^{-4}$ & $1.989624\times10^{-4}$ & $3.289656\times10^{-4}$ & \\
					\midrule
					
					\multirow{3}{*}{1238}
					& $x$ & $5.162852\times10^{-3}$ & $5.186230\times10^{-3}$ & $5.921158\times10^{-3}$ & \multirow{3}{*}{$3.713242\times10^{-3}$} \\
					& $y$ & $1.331224\times10^{-3}$ & $1.351524\times10^{-3}$ & $1.674414\times10^{-3}$ & \\
					& $z$ & $3.552884\times10^{-3}$ & $3.555404\times10^{-3}$ & $3.729424\times10^{-3}$ & \\
					\midrule
					
					\multirow{3}{*}{1239}
					& $x$ & $1.905061\times10^{-4}$ & $2.423644\times10^{-4}$ & $4.548552\times10^{-4}$ & \multirow{3}{*}{$1.499765\times10^{-4}$} \\
					& $y$ & $4.927924\times10^{-5}$ & $5.595751\times10^{-5}$ & $8.532256\times10^{-5}$ & \\
					& $z$ & $5.685171\times10^{-5}$ & $7.488056\times10^{-5}$ & $1.519464\times10^{-4}$ & \\
					\midrule
					
					\multirow{3}{*}{1240}
					& $x$ & $2.048747\times10^{-4}$ & $2.624104\times10^{-4}$ & $4.853309\times10^{-4}$ & \multirow{3}{*}{$1.611155\times10^{-4}$} \\
					& $y$ & $4.898126\times10^{-5}$ & $5.419302\times10^{-5}$ & $8.166822\times10^{-5}$ & \\
					& $z$ & $5.842658\times10^{-5}$ & $7.796503\times10^{-5}$ & $1.593409\times10^{-4}$ & \\
					\midrule
					
					\multirow{3}{*}{1241}
					& $x$ & $2.516110\times10^{-3}$ & $2.576858\times10^{-3}$ & $3.237878\times10^{-3}$ & \multirow{3}{*}{$1.726738\times10^{-3}$} \\
					& $y$ & $5.398933\times10^{-4}$ & $5.412438\times10^{-4}$ & $5.860503\times10^{-4}$ & \\
					& $z$ & $1.397981\times10^{-3}$ & $1.418353\times10^{-3}$ & $1.695860\times10^{-3}$ & \\
					\midrule
					
					\multirow{3}{*}{1242}
					& $x$ & $1.551430\times10^{-4}$ & $1.901201\times10^{-4}$ & $3.664752\times10^{-4}$ & \multirow{3}{*}{$1.207510\times10^{-4}$} \\
					& $y$ & $4.965467\times10^{-5}$ & $5.578220\times10^{-5}$ & $8.464544\times10^{-5}$ & \\
					& $z$ & $5.297984\times10^{-5}$ & $6.697122\times10^{-5}$ & $1.326309\times10^{-4}$ & \\
					\midrule
					
					\multirow{3}{*}{1243}
					& $x$ & $4.394631\times10^{-2}$ & $4.430123\times10^{-2}$ & $4.949229\times10^{-2}$ & \multirow{3}{*}{$3.176514\times10^{-2}$} \\
					& $y$ & $1.305206\times10^{-2}$ & $1.314090\times10^{-2}$ & $1.490195\times10^{-2}$ & \\
					& $z$ & $2.951794\times10^{-2}$ & $2.986286\times10^{-2}$ & $3.465152\times10^{-2}$ & \\
					\bottomrule
				\end{tabular}%
			}
			\vspace{2mm}
			
			\footnotesize
			MAE: mean absolute error; RMSE: root mean squared error; MaxAE: maximum absolute error. The overall RMSE is computed across all three coordinates
			and all sampled parameter values.
		\end{table}
		
		\begin{table}[htbp]
			\centering
			\caption{Comparison at $\tau=2$ (geodesic end points PINN vs Numerical).}
			\label{tab:final_position}
			\begin{tabular}{lccc}
				\toprule
				Method & $x(2)$ & $y(2)$ & $z(2)$ \\
				\midrule
				Numerical & 3.60555128 & 2.58800260 & 3.00000000 \\
				PINN & 3.60543630 & 2.58794780 & 3.00008150 \\
				\bottomrule
			\end{tabular}
		\end{table}
		
		\FloatBarrier
		
		\section{Future Works and Conclusion}
		
		\subsection{Future Works}
		
		Several extensions of this work are possible. The PINN framework may be investigated for more general Riemannian metrics and higher-dimensional manifolds.
		Further improvements in neural network architecture, training strategies, and numerical integration may also help improve the accuracy and stability of the computed geodesics. Future work will additionally focus on systematic grid-refinement studies to establish the convergence of the computed intrinsic area and volume as well as on degenerate metrics. For our experiment we solely focused on positive definite metric, thus experimenting with indefinite metrics, or metrics with singularities will lead us to another dimension of this study.
		
		\subsection{Conclusion}
		
		In this work, we presented a computational framework for visualizing geodesics and geodesic spheres associated with an abstract Riemannian metric on $\mathbb{R}^3$. The geodesic equations were solved using both \texttt{SciPy} and a Physics-Informed Neural Network (PINN), and their solutions were compared to assess
		the accuracy of the neural network approximation. The resulting geodesic sphere was compared with the Euclidean reference sphere using displacement visualization and intrinsic area and volume calculations. It was found that, at the base grid resolution, the computed intrinsic area and volume were approximately $0.362\%$ and $0.364\%$ smaller than their Euclidean reference values, respectively. Grid refinement produced small but nonzero changes in these estimates, indicating the need for further numerical convergence studies.
		
		Overall, this work illustrates how numerical differential equation solvers and PINNs can be combined to study and visualize geodesic structures associated with Riemannian metrics. The approach provides a computational framework for exploring geometric properties that may be difficult to visualize directly.

\vskip6pt
\noindent \textbf{Data availability:} The author(s) declare that all the data generated during the study are in the manuscript and no additional data was generated. The python script used for this study is available in the GitHub repository \cite{python_github}\vskip6pt

\noindent \textbf{Funding:} The author(s) declare that no funding was received for this work.\vskip6pt

\noindent \textbf{Declaration of Generative AI and AI-assisted technologies:} The author(s) declares that no generative AI or AI assisted technology is used for the completion of the manuscript.\vskip6pt

\noindent \textbf{Conflict of interest:} The author(s) declares no conflict of interest.\vskip6pt

\noindent\textbf{Acknowledgement:} 
The author(s) acknowledges all the reviewers for giving their valuable suggestions towards the improvement of the manuscript.


\vspace{0.1in}
\noindent Sujit Bhattacharyya\\
Department of Mathematics, Sister Nivedita University, DG 1/2 New Town, Action Area 1, Kolkata - 700156, India\\
Email: \texttt{sujitbhattacharyya.1996@gmail.com}
	
\end{document}